\documentclass[twoside,11pt]{article}

\usepackage{amsmath,amssymb, amsthm}
\usepackage{graphicx}

\newtheorem{thm}{Theorem}
\newtheorem{cor}{Corollary}
\newtheorem{lem}{Lemma}

\theoremstyle{definition}
\newtheorem{defn}{Definition}
\theoremstyle{remark}
\newtheorem{rem}{Remark}
\numberwithin{equation}{section}
\begin{document}
\pagestyle{myheadings} \markboth{ \rm \centerline {Bogdan Szal}} 
{\rm \centerline { A new class of numerical sequences and its applications to ...}}

\begin{titlepage}
\title{\bf {A new class of numerical sequences and its applications to uniform
convergence of sine series}}
\author {Bogdan Szal \\
{\small University of Zielona G\'{o}ra,}\\
{\small Faculty of Mathematics, Computer Science and Econometrics,}\\
{\small 65-516 Zielona G\'{o}ra, ul. Szafrana 4a, Poland} \\ 
{\small e-mail: B.Szal @wmie.uz.zgora.pl}}
\end{titlepage}

\date{}
\maketitle
\begin{abstract} In the present paper we introduce a new class of sequences called $GM\left(
\beta ,r\right) ,$ which is the generalization of a class considered by
Tikhonov in \cite{12}. Moreover, we obtained in this note sufficient and
necessary conditions for uniform convergence of sine series with $\left(
\beta ,r\right) -$general monotone coefficients.
\end{abstract}



\maketitle

\section{Introduction}

It is well know that there are a great number of interesting results in
Fourier analysis established by assuming monotonicity of coefficients. The
following classical convergence result can be found in many monographs (see 
\cite{1} and \cite{15}, for example).

\begin{thm}
Suppose that $b_{n}\geq b_{n+1}$ and $b_{n}\rightarrow 0$. Then a necessary
and sufficient condition for the uniform convergence of the series%
\begin{equation}
\sum\limits_{n=1}^{\infty }b_{n}\sin nx  \label{1}
\end{equation}%
is $nb_{n}\rightarrow 0$.
\end{thm}

This result have been generalized by weakening the monotone conditions of
the coefficient sequences. Generally speaking, it has become an important
topic how to generalize monotonicity.

Recently, Leindler \cite{3} defined a new class of sequences named as
sequences of rest bounded variation, briefly denoted by $RBVS$, i.e.,%
\begin{equation*}
RBVS=\left\{ a:=\left( a_{n}\right) \in 
\mathbb{C}
:\sum\limits_{n=m}^{\infty }\left\vert a_{n}-a_{n+1}\right\vert \leq
K\left( a\right) \left\vert a_{m}\right\vert \text{ for all }m\in 
\mathbb{N}
\right\} ,
\end{equation*}%
where here and throughout the paper $K\left( a\right) $ always indicates a
constant only depending on $a$.

Denote by $MS$ the class of monotone decreasing sequences and $CQMS$ the
class of classic quasimonotone decreasing sequences ($a\in CQMS$ means that $%
\left( a_{n}\right) \in 
\mathbb{R}
_{+}$ and there exists an $\alpha >0$ such that $a_{n}/n^{\alpha }$ is
decreasing), then it is obvious that%
\begin{equation}
MS\subset RBVS\cap CQMS.  \label{1a}
\end{equation}%
Leindler \cite{16} proved that the class $CQMS$ and $RBVS$ are not
comparable. In \cite{8} leindler considered the class of mean rest bounded
variation sequences $MRBVS$, where%
\begin{equation*}
MRBVS=\left\{ a:=\left( a_{n}\right) \in 
\mathbb{C}
:\right.
\end{equation*}%
\begin{equation*}
\left. \sum\limits_{n=m}^{\infty }\left\vert a_{n}-a_{n+1}\right\vert \leq
K\left( a\right) \frac{1}{m}\sum\limits_{n\geq m/2}^{m}\left\vert
a_{n}\right\vert \text{ for all }m\in 
\mathbb{N}
\right\} .
\end{equation*}%
Further, the class of general monotone coefficients, $GM$, is defined as
follows ( see \cite{10}):%
\begin{equation*}
GM=\left\{ a:=\left( a_{n}\right) \in 
\mathbb{C}
:\sum\limits_{n=m}^{2m-1}\left\vert a_{n}-a_{n+1}\right\vert \leq K\left(
a\right) \left\vert a_{m}\right\vert \text{ for all }m\in 
\mathbb{N}
\right\} .
\end{equation*}%
It is clear that%
\begin{equation}
RBVS\subset MRBVS\text{ \ \ and \ \ }RBVS\cup CQMS\subset GM\text{.}
\label{1b}
\end{equation}

Very recently, Le and Zhou \cite{2} suggested the following new class of
sequences to include $GM$:%
\begin{equation*}
GBVS=\left\{ a:=\left( a_{n}\right) \in 
\mathbb{C}
:\right.
\end{equation*}%
\begin{equation*}
\left. \sum\limits_{n=m}^{2m-1}\left\vert a_{n}-a_{n+1}\right\vert \leq
K\left( a\right) \underset{m\leq n\leq N+m}{\max }\left\vert
a_{n}\right\vert \text{ for some integer }N\text{ and all }m\in 
\mathbb{N}
\right\}
\end{equation*}%
In \cite{6, 10, 11, 12} was defined the class of $\beta -$general monotone
sequences as follows:

\begin{defn}
Let $\beta :=\left( \beta _{n}\right) $ be a nonnegative sequence. The
sequence of complex numbers $a:=\left( a_{n}\right) $ is said to be $\beta -$%
general monotone, or $a\in GM\left( \beta \right) $, if the relation%
\begin{equation*}
\sum\limits_{n=m}^{2m-1}\left\vert a_{n}-a_{n+1}\right\vert \leq K\left(
a\right) \beta _{m}
\end{equation*}%
holds for all $m$.
\end{defn}

In the paper \cite{12} Tikhonov considered the following examples of the
sequences $\beta _{n}:$

(1) $_{1}\beta _{n}=\left\vert a_{n}\right\vert ,$

(2) $_{2}\beta _{n}=\sum\limits_{k=n}^{n+N}\left\vert a_{k}\right\vert $
for some integer $N$,

(3) $_{3}\beta _{n}=\sum\limits_{\nu =0}^{N}\left\vert a_{c^{\nu
}n}\right\vert $ for some integers $N$ and $c>1$,

(4) $_{4}\beta _{n}=\left\vert a_{n}\right\vert +\sum\limits_{k=n+1}^{\left[
cn\right] }\frac{\left\vert a_{k}\right\vert }{k}$ for some $c>1$,

(5) $_{5}\beta _{n}=\left\vert a_{n}\right\vert +\sum\limits_{k=\left[ n/c%
\right] }^{\left[ cn\right] }\frac{\left\vert a_{k}\right\vert }{k}$ for
some $c>1$.

It is clear that $GM\left( _{1}\beta \right) =GM$ and $GM\left( _{2}\beta
\right) =GBVS$. Moreover (see \cite[Remark 2.1]{12})%
\begin{equation}
GM\left( _{1}\beta +_{2}\beta +_{3}\beta +_{4}\beta +_{5}\beta \right)
\equiv GM\left( _{5}\beta \right) .  \label{1c}
\end{equation}

The following results generalized the Chaudy - Joliffe criteria (Theorem 1)
as well as its extensions (see also (\ref{1a}), (\ref{1b}) and (\ref{1c})): 
\cite{9} for $CQMS$, \cite{3} for $RBVS$, \cite{13} for $MRBVS$, \cite{10}
for $GM$, \cite{2} for $GBVS$, \cite{14} for particular case of $GM\left(
_{3}\beta \right) $ sequences, \cite{17, 18} for $GM\left( _{5}\beta \right)
,$ \cite{6} and \cite{11} for $GM\left( \beta
\right) $ but only the sufficient condition and \cite{12} for $GM\left(
\beta \right) $.

We write $I_{1}\ll I_{2}$ if there exists a positive constant $K$ such that $%
I_{1}\leq KI_{2}$.

In order to formulate our new results we define another class of sequences.

\begin{defn}
Let $\beta :=\left( \beta _{n}\right) $ be a nonnegative sequence and $r$ a
natural number. The sequence of complex number $a:=\left( a_{n}\right) $ is
said to be $\left( \beta ,r\right) -$general monotone, or $a\in GM\left(
\beta ,r\right) $, if the relation%
\begin{equation*}
\sum\limits_{n=m}^{2m-1}\left\vert a_{n}-a_{n+r}\right\vert \leq K\left(
a\right) \beta _{m}
\end{equation*}%
holds for all $m$.
\end{defn}

It is clear that $GM\left( \beta ,1\right) \equiv GM\left( \beta \right) $.
Moreover, the embedding relation between $GM\left( \beta ,r\right) $ $\left(
r>1\right) $ and $GM\left( \beta ,1\right) $ implies from the following
remark:

\begin{rem}
Let $r$ be a natural number. If a nonnegative sequence $\beta :=\left( \beta
_{n}\right) $ is such that%
\begin{equation*}
\sum\limits_{i=0}^{r-1}\beta _{n+i}\ll \beta _{n}
\end{equation*}%
for all $n$, then%
\begin{equation*}
GM\left( \beta ,1\right) \subseteq GM\left( \beta ,r\right) .
\end{equation*}
\end{rem}

Connecting with the relation (\ref{1c}) we shall consider the class $%
GM\left( \beta ^{\ast },r\right) $ only, where%
\begin{equation*}
\beta ^{\ast }:=\beta ^{\ast }\left( r\right)
=\sum\limits_{k=n}^{n+r-1}\left\vert a_{k}\right\vert +\sum\limits_{k= 
\left[ n/c\right] }^{\left[ cn\right] }\frac{\left\vert a_{k}\right\vert }{k}%
\text{ for some }c>1.
\end{equation*}%
It is clear that $_{5}\beta \equiv \beta ^{\ast }\left( 1\right) $ and $%
GM\left( _{5}\beta \right) \equiv GM\left( \beta ^{\ast },1\right) \subseteq
GM\left( \beta ^{\ast },r\right) $ for $r\geq 1$ (see Corollary 1).

In this note we shall present the properties of the class $GM\left( \beta
^{\ast },r\right) .$ Moreover, we generalize and extend to the class $%
GM\left( \beta ^{\ast },r\right) $ the results of Tikhonov, which are
included in \cite{12}.

\section{Main results}

We have the following results:

\begin{thm}
Let $r_{1},r_{2}\in 
\mathbb{N}
$ and $r_{1}<r_{2}$. If $r_{1}\mid r_{2}$, then $GM\left( \beta ^{\ast
},r_{1}\right) \varsubsetneq GM\left( \beta ^{\ast },r_{2}\right) .$
\end{thm}

\begin{cor}
If $r\in 
\mathbb{N}
$ and $r>1$, then $GM\left( _{5}\beta \right) =GM\left( \beta ^{\ast
},1\right) \varsubsetneq GM\left( \beta ^{\ast },r\right) $.
\end{cor}

\begin{thm}
Let $r_{1},r_{2}\in 
\mathbb{N}
.$ If $r_{1}\nmid r_{2}$ and $r_{2}\nmid r_{1}$, then the classes $GM\left(
\beta ^{\ast },r_{1}\right) $ and $GM\left( \beta ^{\ast },r_{2}\right) $
are not comparable.
\end{thm}

\begin{thm}
Let a nonnegative sequence $\left( b_{n}\right) \in GM\left( \beta ^{\ast
},r\right) $, where $r\geq 1$. If the series (\ref{1}) converges uniformly
(or if the series (\ref{1}) is the Fourier series of a continuous function),
then $nb_{n}\rightarrow 0$ as $n\rightarrow \infty $.
\end{thm}

\begin{thm}
Let a sequence $\left( b_{n}\right) \in GM\left( \beta ^{\ast },2\right) .$
If $n\left\vert b_{n}\right\vert \rightarrow 0$ as $n\rightarrow \infty $,
then the series (\ref{1}) converges uniformly.
\end{thm}

\begin{cor}
Let a nonnegative sequence $\left( b_{n}\right) \in GM\left( \beta ^{\ast
},2\right) .$ Then the necessary and sufficient condition for series (\ref{1}%
) to be uniformly convergent is $nb_{n}\rightarrow 0$ as $n\rightarrow
\infty $.
\end{cor}

\begin{rem}
If we confine our attention to the class $GM\left( _{5}\beta \right) $ then
by Corollary 1 the Tikhonov result ( see \cite[Remark 2.2.2]{12} ) follows
from Corollary 2.
\end{rem}

\begin{rem}
There exist a real number $x_{0}$ and a sequence $d:=\left( d_{n}\right) \in
GM\left( \beta ^{\ast },3\right) $, with the property $nd_{n}\rightarrow 0$
as $n\rightarrow \infty $, for which the series (\ref{1}) is divergent in $%
x_{0}.$
\end{rem}

\begin{thm}
Let a sequence $\left( b_{n}\right) \in GM\left( \beta ^{\ast },r\right) $,
where $r\geq 3$. If $n\left\vert b_{n}\right\vert \rightarrow 0$ as $%
n\rightarrow \infty $ and 
\begin{equation*}
\sum\limits_{n=1}^{\infty }\sum\limits_{k=1}^{\left[ r/2\right]
}\left\vert b_{r\cdot n+k}-b_{r\cdot n+r-k}\right\vert <\infty ,
\end{equation*}%
then the series (\ref{1}) converges uniformly.
\end{thm}

\begin{rem}
The above result is essential extension of the Tikhonov result (see \cite[%
Theorem 2.1]{12}). Indeed, for any $r\geq 3$ there exists a sequence $%
a:=\left( a_{n}\right) \in GM\left( \beta ^{\ast },r\right) $, with the
properties: $na_{n}\rightarrow 0$ as $n\rightarrow \infty $ and%
\begin{equation}
\sum\limits_{n=1}^{\infty }\sum\limits_{k=1}^{\left[ r/2\right]
}\left\vert a_{r\cdot n+k}-a_{r\cdot n+r-k}\right\vert <\infty ,  \label{r1}
\end{equation}%
which does not belong to the class $GM\left( \beta ^{\ast },2\right) .$
\end{rem}

\section{Lemmas}

\begin{lem}
Let $r\in 
\mathbb{N}
$, $l\in 
\mathbb{Z}
$ and $a:=\left( a_{n}\right) \in 
\mathbb{C}
$. If $x\neq \frac{2l\pi }{r}$, then for all n%
\begin{equation*}
\sum\limits_{k=n}^{2n-1}a_{k}\sin kx=\frac{-1}{2\sin \left( rx/2\right) }%
\left\{ \sum\limits_{k=n}^{2n-1}\left( a_{k}-a_{k+r}\right) \cos \left( k+%
\frac{r}{2}\right) x\right.
\end{equation*}%
\begin{equation}
+\left. \sum\limits_{k=2n}^{2n+r-1}a_{k}\cos \left( k-\frac{r}{2}\right)
x-\sum\limits_{k=n}^{n+r-1}a_{k}\cos \left( k-\frac{r}{2}\right) x\right\} .
\label{l0}
\end{equation}
\end{lem}

\begin{proof}
An elementary calculation gives%
\begin{equation*}
\sum\limits_{k=n}^{2n-1}a_{k}\sin kx=\sum\limits_{k=n}^{2n-1}\left(
a_{k}-a_{k+r}\right) \sin kx+\sum\limits_{k=n}^{2n-1}a_{k+r}\sin kx
\end{equation*}%
\begin{equation}
=\sum\limits_{k=n}^{2n-1}\left( a_{k}-a_{k+r}\right) \sin kx+\cos
rx\sum\limits_{k=n+r}^{2n+r-1}a_{k}\sin kx-\sin
rx\sum\limits_{k=n+r}^{2n+r-1}a_{k}\cos kx.  \label{l1}
\end{equation}%
On the other hand%
\begin{equation*}
\sum\limits_{k=n}^{2n-1}a_{k}\cos kx=\sum\limits_{k=n}^{2n-1}\left(
a_{k}-a_{k+r}\right) \cos kx+\sum\limits_{k=n}^{2n-1}a_{k+r}\cos kx
\end{equation*}%
\begin{equation*}
=\sum\limits_{k=n}^{2n-1}\left( a_{k}-a_{k+r}\right) \cos kx+\cos
rx\sum\limits_{k=n+r}^{2n+r-1}a_{k}\cos kx+\sin
rx\sum\limits_{k=n+r}^{2n+r-1}a_{k}\sin kx.
\end{equation*}%
Hence%
\begin{equation*}
\left( 1-\cos rx\right) \sum\limits_{k=n+r}^{2n+r-1}a_{k}\cos
kx=\sum\limits_{k=n}^{2n-1}\left( a_{k}-a_{k+r}\right) \cos kx+\sin
rx\sum\limits_{k=n+r}^{2n+r-1}a_{k}\sin kx
\end{equation*}%
\begin{equation*}
-\sum\limits_{k=n}^{n+r-1}a_{k}\cos
kx+\sum\limits_{k=2n}^{2n+r-1}a_{k}\cos kx.
\end{equation*}%
Therefore, if $x\neq \frac{2l\pi }{r}$, then%
\begin{equation*}
\sum\limits_{k=n+r}^{2n+r-1}a_{k}\cos kx=\frac{1}{2\sin \left( rx/2\right) }%
\left\{ \sum\limits_{k=n}^{2n-1}\left( a_{k}-a_{k+r}\right) \cos kx\right.
\end{equation*}%
\begin{equation*}
\left. +\sin rx\sum\limits_{k=n+r}^{2n+r-1}a_{k}\sin
kx-\sum\limits_{k=n}^{n+r-1}a_{k}\cos
kx+\sum\limits_{k=2n}^{2n+r-1}a_{k}\cos kx\right\} .
\end{equation*}%
Putting this to (\ref{l1}) we get%
\begin{equation*}
\sum\limits_{k=n}^{2n-1}a_{k}\sin kx=\sum\limits_{k=n}^{2n-1}\left(
a_{k}-a_{k+r}\right) \sin kx+\cos rx\sum\limits_{k=n+r}^{2n+r-1}a_{k}\sin kx
\end{equation*}%
\begin{equation*}
-\frac{\cos \left( rx/2\right) }{\sin \left( rx/2\right) }%
\sum\limits_{k=n}^{2n-1}\left( a_{k}-a_{k+r}\right) \cos kx-2\cos ^{2}\frac{%
rx}{2}\sum\limits_{k=n+r}^{2n+r-1}a_{k}\sin kx
\end{equation*}%
\begin{equation*}
+\frac{\cos \left( rx/2\right) }{\sin \left( rx/2\right) }%
\sum\limits_{k=n}^{n+r-1}a_{k}\cos kx-\frac{\cos \left( rx/2\right) }{\sin
\left( rx/2\right) }\sum\limits_{k=2n}^{2n+r-1}a_{k}\cos kx
\end{equation*}%
\begin{equation*}
=\frac{-1}{\sin \left( rx/2\right) }\sum\limits_{k=n}^{2n-1}\left(
a_{k}-a_{k+r}\right) \cos \left( k+\frac{r}{2}\right)
-\sum\limits_{k=n+r}^{2n+r-1}a_{k}\sin kx
\end{equation*}%
\begin{equation*}
+\frac{\cos \left( rx/2\right) }{\sin \left( rx/2\right) }%
\sum\limits_{k=n}^{n+r-1}a_{k}\cos kx-\frac{\cos \left( rx/2\right) }{\sin
\left( rx/2\right) }\sum\limits_{k=2n}^{2n+r-1}a_{k}\cos kx.
\end{equation*}%
Thus%
\begin{equation*}
2\sum\limits_{k=n}^{2n-1}a_{k}\sin kx=\frac{-1}{\sin \left( rx/2\right) }%
\sum\limits_{k=n}^{2n-1}\left( a_{k}-a_{k+r}\right) \cos \left( k+\frac{r}{2%
}\right)
\end{equation*}%
\begin{equation*}
+\frac{\cos \left( rx/2\right) }{\sin \left( rx/2\right) }%
\sum\limits_{k=n}^{n+r-1}a_{k}\cos kx-\frac{\cos \left( rx/2\right) }{\sin
\left( rx/2\right) }\sum\limits_{k=2n}^{2n+r-1}a_{k}\cos kx
\end{equation*}%
\begin{equation*}
+\sum\limits_{k=n}^{n+r-1}a_{k}\sin
kx-\sum\limits_{k=2n}^{2n+r-1}a_{k}\sin kx
\end{equation*}%
\begin{equation*}
=\frac{-1}{\sin \left( rx/2\right) }\left\{ \sum\limits_{k=n}^{2n-1}\left(
a_{k}-a_{k+r}\right) \cos \left( k+\frac{r}{2}\right) x\right.
\end{equation*}%
\begin{equation*}
+\left. \sum\limits_{k=2n}^{2n+r-1}a_{k}\cos \left( k-\frac{r}{2}\right)
x-\sum\limits_{k=n}^{n+r-1}a_{k}\cos \left( k-\frac{r}{2}\right) x\right\}
\end{equation*}%
and (\ref{l0}) holds.

The proof is complete.
\end{proof}

\begin{defn}
A complex sequence $d:=\left( d_{n}\right) $ is said to be weak monotone if%
\begin{equation*}
n\left\vert d_{n}\right\vert \leq K\left( d\right) \sum\limits_{k=\left[ n/c%
\right] }^{\left[ cn\right] }\left\vert d_{k}\right\vert \text{, \ \ }c>1,
\end{equation*}%
holds for all $n$.
\end{defn}

\begin{lem}
\cite[Theorem 2.2]{12} Let a nonnegative sequence $\left( b_{n}\right) $ be
weak monotone. Then the uniform convergence of series (\ref{1}) (or the fact
that series (\ref{l0}) is the Fourier series of a continuous function)
implies the condition $nb_{n}\rightarrow 0$ as $n\rightarrow \infty $.
\end{lem}

\section{Proofs}

In this section we shall prove our theorems and remarks.

\subsection{Proof of Remark 1}

Let $r\in 
\mathbb{N}
$ and $\left( a_{n}\right) \in GM\left( \beta ,1\right) .$ Then for all $n$%
\begin{equation*}
\sum\limits_{k=n}^{2n-1}\left\vert a_{k}-a_{k+r}\right\vert
=\sum\limits_{k=n}^{2n-1}\left\vert
\sum\limits_{i=0}^{r-1}\left(a_{k+i}-a_{k+i+1}\right) \right\vert
\end{equation*}%
\begin{equation*}
\leq \sum\limits_{k=n}^{2n-1}\sum\limits_{i=0}^{r-1}\left\vert
a_{k+i}-a_{k+i+1}\right\vert
=\sum\limits_{i=0}^{r-1}\sum\limits_{k=n+i}^{2n+i-1}\left\vert
a_{k}-a_{k+1}\right\vert
\ll \sum\limits_{i=0}^{r-1}\beta _{n+i}\ll \beta _{n}
\end{equation*}%
and $\left( a_{n}\right) \in GM\left( \beta ,r\right) $. $\square $

\subsection{Proof of Theorem 2}

If $r_{1}\mid r_{2}$, then exists a natural number $p$ such that $%
r_{2}=p\cdot r_{1}.$ Supposing $\left( a_{n}\right) \in GM\left( \beta
^{\ast },r_{1}\right) $ we have for all $n$%
\begin{equation*}
\sum\limits_{k=n}^{2n-1}\left\vert a_{k}-a_{k+r_{2}}\right\vert
=\sum\limits_{k=n}^{2n-1}\left\vert \sum\limits_{l=0}^{p-1}a_{k+l\cdot
r_{1}}-a_{k+\left( l+1\right) \cdot r_{1}}\right\vert
\end{equation*}%
\begin{equation*}
\leq \sum\limits_{k=n}^{2n-1}\sum\limits_{l=0}^{p-1}\left\vert a_{k+l\cdot
r_{1}}-a_{k+\left( l+1\right) \cdot r_{1}}\right\vert
=\sum\limits_{l=0}^{p-1}\sum\limits_{k=n+l\cdot r_{1}}^{2n+l\cdot
r_{1}-1}\left\vert a_{k}-a_{k+r_{1}}\right\vert
\end{equation*}%
\begin{equation*}
\ll \sum\limits_{l=0}^{p-1}\left( \sum\limits_{k=n+l\cdot r_{1}}^{n+l\cdot
r_{1}+r_{1}-1}\left\vert a_{k}\right\vert +\sum\limits_{k=\left[
c^{-1}\left( n+l\cdot r_{1}\right) \right] }^{\left[ c\left( n+l\cdot
r_{1}\right) \right] }\frac{\left\vert a_{k}\right\vert }{k}\right)
\end{equation*}%
\begin{equation*}
\leq p\left( \sum\limits_{k=n}^{n+p\cdot r_{1}-1}\left\vert
a_{k}\right\vert +\sum\limits_{k=\left[ c^{-1}n\right] }^{\left[ c\left(
n+p\cdot r_{1}\right) \right] }\frac{\left\vert a_{k}\right\vert }{k}\right)
\ll \sum\limits_{k=n}^{n+r_{2}-1}\left\vert a_{k}\right\vert
+\sum\limits_{k=\left[ \left( c\left( 1+r_{2}\right) \right) ^{-1}n\right]
}^{\left[ c\left( 1+r_{2}\right) n\right] }\frac{\left\vert a_{k}\right\vert 
}{k},
\end{equation*}%
whence $\left( a_{n}\right) \in GM\left( \beta ^{\ast },r_{2}\right) $. Thus 
$GM\left( \beta ^{\ast },r_{1}\right) \subseteq GM\left( \beta ^{\ast
},r_{2}\right) $.

Now, we prove that $GM\left( \beta ^{\ast },r_{1}\right) \neq GM\left( \beta
^{\ast },r_{2}\right) $.

Let%
\begin{equation*}
a_{n}=\frac{2+\alpha _{n}}{n^{2}}\text{, \ \ where \ \ \ }\alpha
_{n}=\left\{ 
\begin{array}{c}
-1\text{ \ \ if \ \ }r_{2}\mid n\text{,} \\ 
1\text{ \ \ if \ \ }r_{2}\nmid n.%
\end{array}%
\right.
\end{equation*}%
We show that $\left( a_{n}\right) \in GM\left( \beta ^{\ast },r_{2}\right) $
and $\left( a_{n}\right) \notin GM\left( \beta ^{\ast },r_{1}\right) .$ Let 
\begin{equation*}
A_{r_{2}}:=A\left( r_{2},k,n\right) =\left\{ k:n\leq k<2n\text{ \ \ and \ \ }%
r_{2}\mid k\right\} ,
\end{equation*}%
\begin{equation*}
B_{r_{2}}:=B\left( r_{2},k,n\right) =\left\{ k:n\leq k<2n\text{ \ \ and \ \ }%
r_{2}\nmid k\right\} .
\end{equation*}%
Then for all $n$%
\begin{equation*}
\sum\limits_{k=n}^{2n-1}\left\vert a_{k}-a_{k+r_{2}}\right\vert =\left(
\sum\limits_{k\in A_{r_{2}}}+\sum\limits_{k\in B_{r_{2}}}\right)
\left\vert a_{k}-a_{k+r_{2}}\right\vert
\end{equation*}%
\begin{equation*}
=\sum\limits_{k\in A_{r_{2}}}\left\vert \frac{1}{k^{2}}-\frac{1}{\left(
k+r_{2}\right) ^{2}}\right\vert +\sum\limits_{k\in B_{r_{2}}}\left\vert 
\frac{3}{k^{2}}-\frac{3}{\left( k+r_{2}\right) ^{2}}\right\vert
\end{equation*}%
\begin{equation*}
\leq 3\sum\limits_{k=n}^{2n-1}\frac{2kr_{2}+r_{2}^{2}}{k^{2}\left(
k+r_{2}\right) ^{2}}\leq 6r_{2}\sum\limits_{k=n}^{2n-1}\frac{1}{k^{3}}\leq 
\frac{6r_{2}}{n^{2}}\ll \sum\limits_{k=\left[ n/c\right] }^{\left[ cn\right]
}\frac{1}{k^{3}}\leq \sum\limits_{k=\left[ n/c\right] }^{\left[ cn\right] }%
\frac{a_{k}}{k}
\end{equation*}%
and we have $\left( a_{n}\right) \in GM\left( \beta ^{\ast },r_{2}\right) .$

Since $r_{1}<r_{2}$ we obtain that $r_{2}\nmid r_{1}$ and%
\begin{equation*}
\sum\limits_{k=n}^{2n-1}\left\vert a_{k}-a_{k+r_{1}}\right\vert \geq
\sum\limits_{k\in A_{r_{2}}}\left\vert a_{k}-a_{k+r_{1}}\right\vert
=\sum\limits_{k\in A_{r_{2}}}\left\vert \frac{1}{k^{2}}-\frac{3}{\left(
k+r_{1}\right) ^{2}}\right\vert
\end{equation*}%
\begin{equation*}
=\sum\limits_{k\in A_{r_{2}}}\frac{\left\vert
2k^{2}-2kr_{1}-r_{1}^{2}\right\vert }{k^{2}\left( k+r_{1}\right) ^{2}}.
\end{equation*}%
If $n\geq 5r_{1}$ then $2n^{2}-2nr_{1}-r_{1}^{2}\geq \left( n+r_{1}\right)
^{2}$. Hence for $n\geq 5r_{1}$%
\begin{equation*}
\sum\limits_{k=n}^{2n-1}\left\vert a_{k}-a_{k+r_{1}}\right\vert \geq
\sum\limits_{k\in A_{r_{2}}}\frac{1}{k^{2}}\geq \frac{1}{4nr_{2}}
\end{equation*}%
and since%
\begin{equation*}
\sum\limits_{k=n}^{n+r_{1}-1}a_{k}+\sum\limits_{k=\left[ n/c\right] }^{%
\left[ cn\right] }\frac{a_{k}}{k}\leq \frac{3\left( r_{1}+c^{4}\right) }{%
n^{2}},
\end{equation*}%
the inequality%
\begin{equation*}
\sum\limits_{k=n}^{2n-1}\left\vert a_{k}-a_{k+r_{1}}\right\vert \leq
K\left( a\right) \left\{
\sum\limits_{k=n}^{n+r_{1}-1}a_{k}+\sum\limits_{k= \left[ n/c\right] }^{%
\left[ cn\right] }\frac{a_{k}}{k}\right\}
\end{equation*}%
does not hold, that is, $\left( a_{n}\right) $ does not belong to $GM\left(
\beta ^{\ast },r_{1}\right) .$

This complete the proof. $\square $

\subsection{Proof of Theorem 3}

From the above proof we can see that there exists a sequence $\left(
a_{n}\right) \in GM\left( \beta ^{\ast },r_{2}\right) $ such that $\left(
a_{n}\right) \notin GM\left( \beta ^{\ast },r_{1}\right) .$

Analogously, let%
\begin{equation*}
e_{n}=\frac{2+\gamma _{n}}{n^{2}},\text{ \ \ where \ \ }\gamma _{n}=\left\{ 
\begin{array}{c}
-1\text{ \ \ if \ \ }r_{1}\mid n, \\ 
1\text{ \ \ if \ \ }r_{1}\nmid n.%
\end{array}%
\right.
\end{equation*}%
We show that $\left( e_{n}\right) \in GM\left( \beta ^{\ast },r_{1}\right) $
and $\left( e_{n}\right) \notin GM\left( \beta ^{\ast },r_{2}\right) $. Let%
\begin{equation*}
A_{r_{1}}:=A\left( r_{1},k,n\right) =\left\{ k:n\leq k<2n\text{ \ \ and \ \ }%
r_{1}\mid k\right\} ,
\end{equation*}%
\begin{equation*}
B_{r_{1}}:=B\left( r_{1},k,n\right) =\left\{ k:n\leq k<2n\text{ \ \ and \ \ }%
r_{1}\nmid k\right\} .
\end{equation*}%
Then for all $n$%
\begin{equation*}
\sum\limits_{k=n}^{2n-1}\left\vert e_{k}-e_{k+r_{1}}\right\vert =\left(
\sum\limits_{k\in A_{r_{1}}}+\sum\limits_{k\in B_{r_{1}}}\right)
\left\vert e_{k}-e_{k+r_{1}}\right\vert
\end{equation*}%
\begin{equation*}
=\sum\limits_{k\in A_{r_{1}}}\left\vert \frac{1}{k^{2}}-\frac{1}{\left(
k+r_{1}\right) ^{2}}\right\vert +\sum\limits_{k\in B_{r_{1}}}\left\vert 
\frac{3}{k^{2}}-\frac{3}{\left( k+r_{1}\right) ^{2}}\right\vert
\end{equation*}%
\begin{equation*}
\leq 3\sum\limits_{k=n}^{2n-1}\frac{2kr_{1}+r_{1}^{2}}{k^{2}\left(
k+r_{1}\right) ^{2}}\leq 6r_{1}\sum\limits_{k=n}^{2n-1}\frac{1}{k^{3}}\leq 
\frac{6r_{1}}{n^{2}}\ll \sum\limits_{k=\left[ n/c\right] }^{\left[ cn\right]
}\frac{1}{k^{3}}\leq \sum\limits_{k=\left[ n/c\right] }^{\left[ cn\right] }%
\frac{e_{k}}{k}
\end{equation*}%
and $\left( e_{n}\right) \in GM\left( \beta ^{\ast },r_{1}\right) .$

If $r_{1}\nmid r_{2}$ then for $n\geq 5r_{2}$ we get%
\begin{equation*}
\sum\limits_{k=n}^{2n-1}\left\vert e_{k}-e_{k+r_{2}}\right\vert \geq
\sum\limits_{k\in A_{r_{1}}}\left\vert e_{k}-e_{k+r_{2}}\right\vert
=\sum\limits_{k\in A_{r_{1}}}\left\vert \frac{1}{k^{2}}-\frac{3}{\left(
k+r_{1}\right) ^{2}}\right\vert
\end{equation*}%
\begin{equation*}
=\sum\limits_{k\in A_{r_{1}}}\frac{\left\vert
2k^{2}-2kr_{2}-r_{2}^{2}\right\vert }{k^{2}\left( k+r_{2}\right) ^{2}}\geq
\sum\limits_{k\in A_{r_{1}}}\frac{1}{k^{2}}\geq \frac{1}{4nr_{1}}
\end{equation*}%
and since%
\begin{equation*}
\sum\limits_{k=n}^{n+r_{2}-1}e_{k}+\sum\limits_{k=\left[ n/c\right] }^{%
\left[ cn\right] }\frac{e_{k}}{k}\leq \frac{3\left( r_{2}+c^{4}\right) }{%
n^{2}},
\end{equation*}%
the inequality%
\begin{equation*}
\sum\limits_{k=n}^{2n-1}\left\vert e_{k}-e_{k+r_{1}}\right\vert \leq
K\left( e\right) \left\{
\sum\limits_{k=n}^{n+r_{2}-1}e_{k}+\sum\limits_{k= \left[ n/c\right] }^{%
\left[ cn\right] }\frac{e_{k}}{k}\right\}
\end{equation*}%
does not hold, that is, $\left( e_{n}\right) $ does not belong to $GM\left(
\beta ^{\ast },r_{2}\right) .$

This complete the proof. $\square $

\subsection{Proof of Theorem 4}

Let $\left( b_{n}\right) \in GM\left( \beta ^{\ast },r\right) $, where $%
r\geq 1$. We show that such sequence $\left( b_{n}\right) $ is weak
monotone. If $n\leq r$ then the inequality obviously holds.

Now, let $n>r$. For $j=n+1,n+2,...,2n$ we get%
\begin{equation*}
\sum\limits_{k=n}^{j-1}\left\vert b_{k}-b_{k+r}\right\vert \geq
\sum\limits_{k=n}^{j-1}b_{k}-\sum\limits_{k=n+r}^{j+r-1}b_{k}
\end{equation*}%
and for $j\geq n+r+1$ we obtain%
\begin{equation*}
b_{n}\leq \sum\limits_{k=n}^{n+r-1}b_{k}\leq
\sum\limits_{k=n}^{j-1}\left\vert b_{k}-b_{k+r}\right\vert
+\sum\limits_{k=j}^{j+r-1}b_{k}
\end{equation*}%
\begin{equation*}
\leq \sum\limits_{k=\left[ j/2\right] }^{2\left[ j/2\right] -1}\left\vert
b_{k}-b_{k+r}\right\vert +\sum\limits_{k=j}^{j+r-1}b_{k}
\end{equation*}%
\begin{equation*}
\ll \sum\limits_{k=\left[ j/2c\right] }^{\left[ cj/2\right] }\frac{b_{k}}{k}%
+\sum\limits_{k=\left[ j/2\right] }^{\left[ j/2\right] +r-1}b_{k}+\sum%
\limits_{k=j}^{j+r-1}b_{k}\leq \frac{1}{n}\sum\limits_{k=\left[ n/2c\right]
}^{\left[ cn\right] }b_{k}+\sum\limits_{k=\left[ j/2\right] }^{\left[ j/2%
\right] +r-1}b_{k}+\sum\limits_{k=j}^{j+r-1}b_{k}.
\end{equation*}%
Summing up on $j$ we get%
\begin{equation*}
nb_{n}=\sum\limits_{j=n+1}^{2n}b_{n}=\sum\limits_{j=n+1}^{n+r}b_{n}+\sum%
\limits_{j=n+r+1}^{2n}b_{n}
\end{equation*}%
\begin{equation*}
\ll rb_{n}+\sum\limits_{j=n+r+1}^{2n}\left( \frac{1}{n}\sum\limits_{k=%
\left[ n/2c\right] }^{\left[ cn\right] }b_{k}+\sum\limits_{k=\left[ j/2%
\right] }^{\left[ j/2\right] +r-1}b_{k}+\sum\limits_{k=j}^{j+r-1}b_{k}%
\right)
\end{equation*}%
\begin{equation*}
\leq rb_{n}+\sum\limits_{k=\left[ n/2c\right] }^{\left[ cn\right]
}b_{k}+\sum\limits_{j=n+1}^{2n}\sum\limits_{k=j}^{j+r-1}b_{k}+\sum%
\limits_{j=n+1}^{2n}\sum\limits_{k=\left[ j/2\right] }^{\left[ j/2\right]
+r-1}b_{k}
\end{equation*}%
\begin{equation*}
\leq r\sum\limits_{k=\left[ n/2c\right] }^{\left[ cn\right]
}b_{k}+\sum\limits_{j=0}^{r-1}\sum\limits_{k=n+1+j}^{2n-1+j}b_{k}+2\sum%
\limits_{j=0}^{r-1}\sum\limits_{k=\left[ \left( n+1\right) /2\right] +j}^{%
\left[ \left( 2n-1\right) /2\right] +j}b_{k}
\end{equation*}%
\begin{equation*}
\leq r\sum\limits_{k=\left[ n/2c\right] }^{\left[ cn\right]
}b_{k}+r\sum\limits_{k=n+1}^{2n+r-2}b_{k}+2r\sum\limits_{k=\left[ \left(
n+1\right) /2\right] }^{\left[ \left( 2n-1\right) /2\right] +r-1}b_{k}
\end{equation*}%
\begin{equation*}
\leq r\sum\limits_{k=\left[ n/2c\right] }^{\left[ cn\right]
}b_{k}+3r\sum\limits_{k=\left[ n/2\right] }^{2n+r}b_{k}\leq r\sum\limits_{k=%
\left[ n/2c\right] }^{\left[ 2cn\right] }b_{k}+3r\sum\limits_{k=\left[ n/2c%
\right] }^{3n}b_{k}\leq 4r\sum\limits_{k=\left[ n/c_{1}\right] }^{\left[ nc_{1}%
\right] }b_{k},
\end{equation*}%
where $c_{1}=\max \left\{ 3,2c\right\} .$ Therefore $\left( b_{n}\right) $
is weak monotone and by Lemma 2 we have that $nb_{n}\rightarrow 0$ as $%
n\rightarrow \infty $ and thus the proof is complete. $\square $

\subsection{Proof of Theorem 5}

Denote%
\begin{equation*}
\varepsilon _{n}^{\left( 1\right) }:=\underset{k\geq n/c}{\sup }k\left\vert
b_{k}\right\vert \text{, \ \ where \ \ }c>1
\end{equation*}%
and%
\begin{equation*}
r_{n}\left( x\right) =\sum\limits_{k=n}^{\infty }b_{k}\sin kx.
\end{equation*}%
In view of the assumption, we can see that $\varepsilon _{n}^{\left(
1\right) }\rightarrow 0$ as $n\rightarrow \infty $. Further, we shall show
that 
\begin{equation}
\left\vert r_{n}\left( x\right) \right\vert \ll \varepsilon _{n}^{\left(
1\right) }.  \label{t1}
\end{equation}%
Since $r_{n}\left( l\pi \right) =0$, where $l\in 
\mathbb{Z}
$, it suffices to prove (\ref{t1}) for $0<x<\pi $.

First we show that (\ref{t1}) is valid for $0<x\leq \frac{\pi }{2}$.

Let $N:=N\left( x\right) \geq 2$ be the natural number such that 
\begin{equation}
\frac{\pi }{N+1}<x\leq \frac{\pi }{N}.  \label{t2}
\end{equation}%
Then%
\begin{equation*}
r_{n}\left( x\right) =\sum\limits_{k=n}^{n+N-1}b_{k}\sin
kx+\sum\limits_{k=n+N}^{\infty }b_{k}\sin kx=r_{n}^{\left( 1\right) }\left(
x\right) +r_{n}^{\left( 2\right) }\left( x\right) .
\end{equation*}%
Hence, by (\ref{t2}),%
\begin{equation}
\left\vert r_{n}^{\left( 1\right) }\left( x\right) \right\vert \leq
\sum\limits_{k=n}^{n+N-1}\left\vert b_{k}\right\vert \left\vert \sin
kx\right\vert \leq x\sum\limits_{k=n}^{n+N-1}k\left\vert b_{k}\right\vert
\leq x\varepsilon _{n}^{\left( 1\right) }N\leq \pi \varepsilon _{n}^{\left(
1\right) }.  \label{t4}
\end{equation}%
If $\left( b_{n}\right) \in GM\left( \beta ^{\ast },2\right) $, then using
Lemma 1, the inequality $\frac{2}{\pi }x\leq \sin x$ $\left( x\in \left[ 0,%
\frac{\pi }{2}\right] \right) $ and (\ref{t2}) we obtain%
\begin{equation*}
\left\vert r_{n}^{\left( 2\right) }\left( x\right) \right\vert =\left\vert
\sum\limits_{j=0}^{\infty }\sum\limits_{k=2^{j}\left( n+N\right)
}^{2^{j+1}\left( n+N\right) -1}b_{k}\sin kx\right\vert
\end{equation*}%
\begin{equation*}
\leq \sum\limits_{j=0}^{\infty }\left\vert \frac{-1}{2\sin x}%
\sum\limits_{k=2^{j}\left( n+N\right) }^{2^{j+1}\left( n+N\right) -1}\left(
b_{k}-b_{k+2}\right) \cos \left( k+1\right) x\right.
\end{equation*}%
\begin{equation*}
\left. \left. +\sum\limits_{k=2^{j+1}\left( n+N\right) }^{2^{j+1}\left(
n+N\right) +1}b_{k}\cos \left( k-1\right) x-\sum\limits_{k=2^{j}\left(
n+N\right) }^{2^{j}\left( n+N\right) +1}b_{k}\cos \left( k-1\right)
x\right\} \right\vert
\end{equation*}%
\begin{equation*}
\leq \frac{1}{2\sin x}\sum\limits_{j=0}^{\infty }\left\{
\sum\limits_{k=2^{j}\left( n+N\right) }^{2^{j+1}\left( n+N\right)
-1}\left\vert b_{k}-b_{k+2}\right\vert +\sum\limits_{k=2^{j+1}\left(
n+N\right) }^{2^{j+1}\left( n+N\right) +1}\left\vert b_{k}\right\vert
+\sum\limits_{k=2^{j}\left( n+N\right) }^{2^{j}\left( n+N\right)
+1}\left\vert b_{k}\right\vert \right\}
\end{equation*}%
\begin{equation*}
\ll \frac{\pi }{2x}\sum\limits_{j=0}^{\infty }\left\{
\sum\limits_{k=2^{j}\left( n+N\right) }^{2^{j+1}\left( n+N\right)
-1}\left\vert b_{k}-b_{k+2}\right\vert +\sum\limits_{k=2^{j}\left(
n+N\right) }^{2^{j}\left( n+N\right) +1}\left\vert b_{k}\right\vert \right\}
\end{equation*}%
\begin{equation*}
\ll \left( N+1\right) \sum\limits_{j=0}^{\infty }\left\{
\sum\limits_{k=2^{j}\left( n+N\right) }^{2^{j}\left( n+N\right)
+1}\left\vert b_{k}\right\vert +\sum\limits_{k=\left[ 2^{j}\left(
n+N\right) /c\right] }^{\left[ c2^{j}\left( n+N\right) \right] }\frac{%
\left\vert b_{k}\right\vert }{k}\right\}
\end{equation*}%
\begin{equation*}
\ll \left( N+1\right) \varepsilon _{n}^{\left( 1\right)
}\sum\limits_{j=0}^{\infty }\left\{ \sum\limits_{k=2^{j}\left( n+N\right)
}^{2^{j}\left( n+N\right) +1}\frac{1}{k}+\sum\limits_{k=\left[ 2^{j}\left(
n+N\right) /c\right] }^{\left[ c2^{j}\left( n+N\right) \right] }\frac{1}{%
k^{2}}\right\}
\end{equation*}%
\begin{equation}
\ll \frac{N+1}{N+n}\varepsilon _{n}^{\left( 1\right)
}\sum\limits_{j=0}^{\infty }\frac{1}{2^{j}}\ll \varepsilon _{n}^{\left(
1\right) }.  \label{t5}
\end{equation}

Now, we prove (\ref{t1}) for $\frac{\pi }{2}\leq x<\pi $.

Let $M:=M\left( x\right) \geq 2$ be the natural number such that%
\begin{equation}
\pi -\frac{\pi }{M}\leq x<\pi -\frac{\pi }{M+1}.  \label{t3}
\end{equation}%
Then%
\begin{equation*}
r_{n}\left( x\right) =\sum\limits_{k=n}^{n+M-1}b_{k}\sin
kx+\sum\limits_{k=n+M}^{\infty }b_{k}\sin kx=r_{n}^{\left( 3\right) }\left(
x\right) +r_{n}^{\left( 4\right) }\left( x\right) .
\end{equation*}%
Using the inequality $\sin x\leq \pi -x$ $\left( x\in \left( 0,\pi \right)
\right) $ and (\ref{t3}) we get%
\begin{equation}
\left\vert r_{n}^{\left( 3\right) }\left( x\right) \right\vert \leq \left(
\pi -x\right) \sum\limits_{k=n}^{n+M-1}k\left\vert b_{k}\right\vert \leq
\left( \pi -x\right) M\varepsilon _{n}^{\left( 1\right) }\leq \pi
\varepsilon _{n}^{\left( 1\right) }.  \label{t6}
\end{equation}%
If $\left( b_{n}\right) \in GM\left( \beta ^{\ast },2\right) $, then using
Lemma 1, the inequality $2-\frac{2}{\pi }x\leq \sin x$ $\left( x\in \left[ 
\frac{\pi }{2},\pi \right] \right) $ and (\ref{t3}) we obtain%
\begin{equation*}
\left\vert r_{n}^{\left( 2\right) }\left( x\right) \right\vert =\left\vert
\sum\limits_{j=0}^{\infty }\sum\limits_{k=2^{j}\left( n+M\right)
}^{2^{j+1}\left( n+M\right) -1}b_{k}\sin kx\right\vert
\end{equation*}%
\begin{equation*}
\leq \sum\limits_{j=0}^{\infty }\left\vert \frac{-1}{2\sin x}%
\sum\limits_{k=2^{j}\left( n+M\right) }^{2^{j+1}\left( n+M\right) -1}\left(
b_{k}-b_{k+2}\right) \cos \left( k+1\right) x\right.
\end{equation*}%
\begin{equation*}
\left. \left. +\sum\limits_{k=2^{j+1}\left( n+M\right) }^{2^{j+1}\left(
n+M\right) +1}b_{k}\cos \left( k-1\right) x-\sum\limits_{k=2^{j}\left(
n+M\right) }^{2^{j}\left( n+M\right) +1}b_{k}\cos \left( k-1\right)
x\right\} \right\vert
\end{equation*}%
\begin{equation*}
\ll \frac{1}{2\left( 1-\frac{1}{\pi }x\right) }\sum\limits_{j=0}^{\infty
}\left\{ \sum\limits_{k=2^{j}\left( n+M\right) }^{2^{j+1}\left( n+M\right)
-1}\left\vert b_{k}-b_{k+2}\right\vert +\sum\limits_{k=2^{j}\left(
n+M\right) }^{2^{j}\left( n+M\right) +1}\left\vert b_{k}\right\vert \right\}
\end{equation*}%
\begin{equation*}
\ll \left( M+1\right) \sum\limits_{j=0}^{\infty }\left\{
\sum\limits_{k=2^{j}\left( n+M\right) }^{2^{j}\left( n+M\right)
+1}\left\vert b_{k}\right\vert +\sum\limits_{k=\left[ 2^{j}\left(
n+M\right) /c\right] }^{\left[ c2^{j}\left( n+M\right) \right] }\frac{%
\left\vert b_{k}\right\vert }{k}\right\}
\end{equation*}%
\begin{equation}
\ll \frac{M+1}{M+n}\varepsilon _{n}^{\left( 1\right)
}\sum\limits_{j=0}^{\infty }\frac{1}{2^{j}}\ll \varepsilon _{n}^{\left(
1\right) }.  \label{t7}
\end{equation}%
From the estimations (\ref{t4}), (\ref{t5}), (\ref{t6}) and (\ref{t7})we
obtain the uniform convergence of series (\ref{1}) follows and thus the
proof is complete. $\square $

\subsection{Proof of Remark 3}

Let $x_{0}=\frac{2\pi }{3}$ and%
\begin{equation*}
d_{n}=\left\{ 
\begin{array}{c}
\frac{3}{n\ln \left( n+1\right) }\text{ \ \ if \ \ }n=3k+1, \\ 
\frac{1}{n\ln \left( n+1\right) }\text{ \ \ if \ \ }n\neq 3k+1,%
\end{array}%
\right.
\end{equation*}%
where $k\in 
\mathbb{N}
\cup \left\{ 0\right\} $. An elementary calculation gives%
\begin{equation*}
\sum\limits_{k=1}^{\infty }d_{k}\sin \left( kx_{0}\right) =d_{1}\sin \frac{%
2\pi }{2}+d_{2}\sin \frac{4\pi }{2}+\sum\limits_{k=1}^{\infty
}\sum\limits_{i=0}^{2}d_{3k+i}\sin \left( 3k+i\right) \frac{2\pi }{3}
\end{equation*}%
\begin{equation*}
=\left( d_{1}-d_{2}\right) \sin \frac{2\pi }{3}+\sum\limits_{k=1}^{\infty
}\sum\limits_{i=0}^{2}d_{3k+i}\sin \left( i\frac{2\pi }{3}\right)
\end{equation*}%
\begin{equation*}
=\left( d_{1}-d_{2}\right) \sin \frac{2\pi }{3}+\sum\limits_{k=1}^{\infty
}\sum\limits_{i=1}^{2}d_{3k+i}\sin \left( i\frac{2\pi }{3}\right) =\sin 
\frac{2\pi }{3}\sum\limits_{k=1}^{\infty }\left( d_{3k+1}-d_{3k+2}\right)
\end{equation*}%
\begin{equation*}
=\sin \frac{2\pi }{3}\sum\limits_{k=1}^{\infty }\left( \frac{3}{\left(
3k+1\right) \ln \left( 3k+2\right) }-\frac{1}{\left( 3k+2\right) \ln \left(
3k+3\right) }\right)
\end{equation*}%
\begin{equation*}
\geq \sin \frac{2\pi }{3}\sum\limits_{k=1}^{\infty }\frac{2}{\left(
3k+1\right) \ln \left( 3k+3\right) }=\infty
\end{equation*}%
and the series $\sum\limits_{k=1}^{\infty }d_{k}\sin \left( kx_{0}\right) $
is divergent, too. $\square $

\subsection{Proof of Theorem 6}

Denote%
\begin{equation*}
\varepsilon _{n}^{\left( 1\right) }:=\underset{k\geq n/c}{\sup }k\left\vert
b_{k}\right\vert \text{, \ \ where \ \ }c>1,
\end{equation*}%
\begin{equation*}
\varepsilon _{n}^{\left( 2\right) }:=\sum\limits_{k\geq n/r}^{\infty
}\sum\limits_{k=1}^{\left[ r/2\right] }\left\vert b_{r\cdot n+k}-b_{r\cdot
n+r-k}\right\vert
\end{equation*}%
and%
\begin{equation*}
r_{n}\left( x\right) =\sum\limits_{k=n}^{\infty }b_{k}\sin kx.
\end{equation*}%
In view of the assumptions, we have that $\varepsilon _{n}^{\left( 1\right)
}\rightarrow 0$ and $\varepsilon _{n}^{\left( 2\right) }\rightarrow 0$ as $%
n\rightarrow \infty $. We shall show that 
\begin{equation}
\left\vert r_{n}\left( x\right) \right\vert \ll \varepsilon _{n}^{\left(
1\right) }+\varepsilon _{n}^{\left( 2\right) }  \label{p1}
\end{equation}%
also holds. Since $r_{n}\left( l\pi \right) =0$, where $l\in 
\mathbb{Z}
$, it suffices to prove (\ref{p1}) for $0<x<\pi $.

First we shall that (\ref{p1}) is valid for $x=\frac{2l\pi }{r}$, where $l$
is an integer number such that $0<2l<r.$ For any $n$ exist two number $%
p,q\in 
\mathbb{N}
\cup \left\{ 0\right\} $ such that $n=r\cdot p+q$, where $0\leq q<r$. Then%
\begin{equation*}
r_{n}\left( \frac{2l\pi }{r}\right) =\sum\limits_{k=n}^{\infty }b_{k}\sin
\left( k\frac{2l\pi }{r}\right) =\sum\limits_{k=r\cdot p+q}^{r\left(
p+1\right) -1}b_{k}\sin \left( k\frac{2l\pi }{r}\right)
\end{equation*}%
\begin{equation*}
+\sum\limits_{k=p+1}^{\infty }\sum\limits_{i=0}^{r-1}b_{r\cdot k+i}\sin
\left( \left( rk+i\right) \frac{2l\pi }{r}\right) =\sum\limits_{k=r\cdot
p+q}^{r\left( p+1\right) -1}b_{k}\sin \left( k\frac{2l\pi }{r}\right)
\end{equation*}%
\begin{equation*}
+\sum\limits_{k=p+1}^{\infty }\sum\limits_{i=1}^{r-1}b_{r\cdot k+i}\sin
\left( i\frac{2l\pi }{r}\right) .
\end{equation*}%
If $r=2s$ $\left( s=2,3,..\right) $ then%
\begin{equation*}
\sum\limits_{i=1}^{r-1}b_{r\cdot k+i}\sin \left( i\frac{2l\pi }{r}\right)
=\sum\limits_{k=1}^{2s-1}b_{2s\cdot k+i}\sin \left( i\frac{l\pi }{s}\right)
=\sum\limits_{i=1}^{s-1}\left( b_{2s\cdot k+i}-b_{2s\cdot k+2s-i}\right)
\sin \left( i\frac{l\pi }{s}\right)
\end{equation*}%
\begin{equation}
=\sum\limits_{i=1}^{s}\left( b_{2s\cdot k+i}-b_{2s\cdot k+2s-i}\right) \sin
\left( i\frac{l\pi }{s}\right) =\sum\limits_{i=1}^{r/2}\left( b_{r\cdot
k+i}-b_{r\cdot k+r-i}\right) \left( i\frac{2l\pi }{r}\right)  \label{p2}
\end{equation}%
and if $r=2s+1$ $\left( s=1,2,..\right) $ then%
\begin{equation*}
\sum\limits_{i=1}^{r-1}b_{r\cdot k+i}\sin \left( i\frac{2l\pi }{r}\right)
=\sum\limits_{k=1}^{2s}b_{\left( 2s+1\right) \cdot k+i}\sin \left( i\frac{%
2l\pi }{2s+1}\right)
\end{equation*}%
\begin{equation*}
=\sum\limits_{k=1}^{s}\left( b_{\left( 2s+1\right) \cdot k+i}-b_{\left(
2s+1\right) \cdot k+2s+1-i}\right) \left( i\frac{2l\pi }{2s+1}\right)
\end{equation*}%
\begin{equation}
=\sum\limits_{i=1}^{\left[ r/2\right] }\left( b_{r\cdot k+i}-b_{r\cdot
k+r-i}\right) \left( i\frac{2l\pi }{r}\right) .  \label{p3}
\end{equation}%
Using (\ref{p2}) or (\ref{p3}) we obtain%
\begin{equation*}
\left\vert r_{n}\left( \frac{2l\pi }{r}\right) \right\vert \leq
\sum\limits_{k=r\cdot p+q}^{r\left( p+1\right) -1}\left\vert
b_{k}\right\vert +\sum\limits_{k=p+1}^{\infty }\sum\limits_{i=1}^{\left[
r/2\right] }\left\vert b_{r\cdot k+i}-b_{r\cdot k+r-i}\right\vert
\end{equation*}%
\begin{equation}
\leq \frac{r}{n}\varepsilon _{n}^{\left( 1\right) }+\varepsilon _{n}^{\left(
2\right) }\ll \varepsilon _{n}^{\left( 1\right) }+\varepsilon _{n}^{\left(
2\right) }.  \label{p9}
\end{equation}

Now, we prove that (\ref{p1}) holds for $\frac{2l\pi }{r}<x\leq \frac{2l\pi 
}{r}+\frac{\pi }{r}$, where $0\leq 2l<r$.

Let $N:=N\left( x\right) $ be the natural number such that%
\begin{equation}
\frac{2l\pi }{r}+\frac{\pi }{N+1}<x\leq \frac{2l\pi }{r}+\frac{\pi }{N}.
\label{p4}
\end{equation}%
Then%
\begin{equation*}
r_{n}\left( x\right) =\sum\limits_{k=n}^{n+N-1}b_{k}\sin
kx+\sum\limits_{k=n+N}^{\infty }b_{k}\sin kx=r_{n}^{\left( 1\right) }\left(
x\right) +r_{n}^{\left( 2\right) }\left( x\right) .
\end{equation*}%
Applying Lagrange's mean value theorem to the function $f\left( x\right)
=\sin kx$ on the interval $\left[ \frac{2l\pi }{r},x\right] $ we obtain that
there exists $y\in \left( \frac{2l\pi }{r},x\right) $ such that%
\begin{equation*}
\sin kx-\sin \left( k\frac{2l\pi }{r}\right) =k\cos ky\left( x-\frac{2l\pi }{%
r}\right) .
\end{equation*}%
Using this we get%
\begin{equation*}
r_{n}^{\left( 1\right) }\left( x\right)
=\sum\limits_{k=n}^{n+N-1}kb_{k}\cos ky\left( x-\frac{2l\pi }{r}\right)
+\sum\limits_{k=n}^{n+N-1}b_{k}\sin \left( k\frac{2l\pi }{r}\right)
\end{equation*}%
\begin{equation*}
=r_{n}^{\left( 1.1\right) }\left( x\right) +r_{n}^{\left( 1.2\right) }\left( 
\frac{2l\pi }{r}\right) .
\end{equation*}%
Hence, by (\ref{p4}),%
\begin{equation}
\left\vert r_{n}^{\left( 1.1\right) }\left( x\right) \right\vert \leq \left(
x-\frac{2l\pi }{r}\right) \sum\limits_{k=n}^{n+N-1}k\left\vert
b_{k}\right\vert \leq \pi \varepsilon _{n}^{\left( 1\right) }.  \label{p10}
\end{equation}%
If $l=0$ then $r_{n}^{\left( 1.2\right) }\left( 0\right) =0$ and (\ref{p1})
is evident. Let $0<2l<r.$ For any $n$ and $N\geq r$ there exist four numbers 
$p_{1},p_{2},q_{1},q_{2}\in 
\mathbb{N}
\cup \left\{ 0\right\} $ such that $n=rp_{1}+q_{1}$ and $N=rp_{2}+q_{2}$,
where $p_{2}\geq 1$ and $0\leq q_{1},q_{2}<r$. If $q_{1}+q_{2}=r$ then%
\begin{equation*}
r_{n}^{\left( 1.2\right) }\left( \frac{2l\pi }{r}\right)
=\sum\limits_{k=r\cdot p_{1+q_{1}}}^{r\left( p_{1}+1\right) -1}b_{k}\sin
\left( k\frac{2l\pi }{r}\right) +\sum\limits_{k=r\left( p_{1}+1\right)
}^{r\left( p_{1}+p_{2}\right) +r-1}b_{k}\sin \left( k\frac{2l\pi }{r}\right)
\end{equation*}%
\begin{equation}
=\sum\limits_{k=r\cdot p_{1+q_{1}}}^{r\left( p_{1}+1\right) -1}b_{k}\sin
\left( k\frac{2l\pi }{r}\right)
+\sum\limits_{k=p_{1}+1}^{p_{1}+p_{2}}\sum\limits_{i=0}^{r-1}b_{r\cdot
k+i}\sin \left( i\frac{2l\pi }{r}\right)  \label{p5}
\end{equation}%
and using (\ref{p2}) or (\ref{p3}) we obtain%
\begin{equation*}
\left\vert r_{n}^{\left( 1.2\right) }\left( \frac{2l\pi }{r}\right)
\right\vert \leq \sum\limits_{k=r\cdot p_{1}+q_{1}}^{r\left( p_{1}+1\right)
-1}\left\vert b_{k}\right\vert
+\sum\limits_{k=p_{1}+1}^{p_{1}+p_{2}}\sum\limits_{i=1}^{\left[ r/2\right]
}\left\vert b_{r\cdot k+i}-b_{r\cdot k+r-i}\right\vert \ll \varepsilon
_{n}^{\left( 1\right) }+\varepsilon _{n}^{\left( 2\right) }.
\end{equation*}%
Let $q_{1}+q_{2}<r.$ Then%
\begin{equation*}
r_{n}^{\left( 1.2\right) }\left( \frac{2l\pi }{r}\right)
=\sum\limits_{k=r\cdot p_{1}+q_{1}}^{r\left( p_{1}+1\right) -1}b_{k}\sin
\left( k\frac{2l\pi }{r}\right)
\end{equation*}%
\begin{equation*}
+\sum\limits_{k=p_{1}+1}^{p_{1}+p_{2}}\sum\limits_{i=0}^{r-1}b_{r\cdot
k+i}\sin \left( i\frac{2l\pi }{r}\right) -\sum\limits_{k=r\left(
p_{1}+p_{2}\right) +q_{1}+q_{2}-1}^{r\left( p_{1}+p_{2}\right)
+r-1}b_{k}\sin \left( k\frac{2l\pi }{r}\right)
\end{equation*}%
and by (\ref{p2}) or (\ref{p3})%
\begin{equation*}
\left\vert r_{n}^{\left( 1.2\right) }\left( \frac{2l\pi }{r}\right)
\right\vert \leq \sum\limits_{k=r\cdot p+q_{1}}^{r\left( p_{1}+1\right)
-1}\left\vert b_{k}\right\vert
+\sum\limits_{k=p_{1}+1}^{p_{1}+p_{2}}\sum\limits_{i=1}^{\left[ r/2\right]
}\left\vert b_{r\cdot k+i}-b_{r\cdot k+r-i}\right\vert
\end{equation*}%
\begin{equation}
+\sum\limits_{k=r\left( p_{1}+p_{2}\right) +q_{1}+q_{2}-1}^{r\left(
p_{1}+p_{2}\right) +r-1}\left\vert b_{k}\right\vert \ll \varepsilon
_{n}^{\left( 1\right) }+\varepsilon _{n}^{\left( 2\right) }.  \label{p6}
\end{equation}%
Now, assume $q_{1}+q_{2}>r$. Then%
\begin{equation*}
r_{n}^{\left( 1.2\right) }\left( \frac{2l\pi }{r}\right)
=\sum\limits_{k=r\cdot p_{1+q_{1}}}^{r\left( p_{1}+1\right) -1}b_{k}\sin
\left( k\frac{2l\pi }{r}\right)
\end{equation*}%
\begin{equation*}
+\sum\limits_{k=p_{1}+1}^{p_{1}+p_{2}}\sum\limits_{i=0}^{r-1}b_{r\cdot
k+i}\sin \left( i\frac{2l\pi }{r}\right) +\sum\limits_{k=r\left(
p_{1}+p_{2}\right) +r}^{r\left( p_{1}+p_{2}\right) +q_{1}+q_{2}-1}b_{k}\sin
\left( k\frac{2l\pi }{r}\right)
\end{equation*}%
and using (\ref{p2}) or (\ref{p3}) we have%
\begin{equation*}
\left\vert r_{n}^{\left( 1.2\right) }\left( \frac{2l\pi }{r}\right)
\right\vert \leq \sum\limits_{k=r\cdot p_{1}+q_{1}}^{r\left( p_{1}+1\right)
-1}\left\vert b_{k}\right\vert
+\sum\limits_{k=p_{1}+1}^{p_{1}+p_{2}}\sum\limits_{i=1}^{\left[ r/2\right]
}\left\vert b_{r\cdot k+i}-b_{r\cdot k+r-i}\right\vert
\end{equation*}%
\begin{equation}
+\sum\limits_{k=r\left( p_{1}+p_{2}\right) +r}^{r\left( p_{1}+p_{2}\right)
+q_{1}+q_{2}-1}\left\vert b_{k}\right\vert \ll \varepsilon _{n}^{\left(
1\right) }+\varepsilon _{n}^{\left( 2\right) }.  \label{p7}
\end{equation}%
Therefore, by (\ref{p5})-(\ref{p6}),%
\begin{equation}
\left\vert r_{n}^{\left( 1.2\right) }\left( \frac{2l\pi }{r}\right)
\right\vert \ll \varepsilon _{n}^{\left( 1\right) }+\varepsilon _{n}^{\left(
2\right) }.  \label{p11}
\end{equation}%
If $\left( b_{n}\right) \in GM\left( \beta ^{\ast },3\right) $ $\left( r\geq
3\right) $, then using Lemma 1, the inequality $\frac{r}{\pi }x-2l\leq
\left\vert \sin \frac{rx}{2}\right\vert $ $\left( x\in \left[ \frac{2l\pi }{r%
},\frac{2l\pi }{r}+\frac{\pi }{r}\right] \text{ and }0\leq 2l<r\right) $ and
(\ref{p4}) we obtain%
\begin{equation*}
\left\vert r_{n}^{\left( 2\right) }\left( x\right) \right\vert =\left\vert
\sum\limits_{j=0}^{\infty }\sum\limits_{k=2^{j}\left( n+N\right)
}^{2^{j+1}\left( n+N\right) -1}b_{k}\sin kx\right\vert
\end{equation*}%
\begin{equation*}
\leq \sum\limits_{j=0}^{\infty }\left\vert \frac{-1}{2\sin \left(
rx/2\right) }\sum\limits_{k=2^{j}\left( n+N\right) }^{2^{j+1}\left(
n+N\right) -1}\left( b_{k}-b_{k+r}\right) \cos \left( k+\frac{r}{2}\right)
x\right.
\end{equation*}%
\begin{equation*}
\left. \left. +\sum\limits_{k=2^{j+1}\left( n+N\right) }^{2^{j+1}\left(
n+N\right) +r-1}b_{k}\cos \left( k-\frac{r}{2}\right)
x-\sum\limits_{k=2^{j}\left( n+N\right) }^{2^{j}\left( n+N\right)
+r-1}b_{k}\cos \left( k-\frac{r}{2}\right) x\right\} \right\vert
\end{equation*}%
\begin{equation*}
\leq \frac{1}{2\left\vert \sin \left( rx/2\right) \right\vert }%
\sum\limits_{j=0}^{\infty }\left\{ \sum\limits_{k=2^{j}\left( n+N\right)
}^{2^{j+1}\left( n+N\right) -1}\left\vert b_{k}-b_{k+r}\right\vert
+\sum\limits_{k=2^{j+1}\left( n+N\right) }^{2^{j+1}\left( n+N\right)
+r-1}\left\vert b_{k}\right\vert +\sum\limits_{k=2^{j}\left( n+N\right)
}^{2^{j}\left( n+N\right) +r-1}\left\vert b_{k}\right\vert \right\}
\end{equation*}%
\begin{equation*}
\ll \frac{1}{\frac{r}{\pi }x-2l}\sum\limits_{j=0}^{\infty }\left\{
\sum\limits_{k=2^{j}\left( n+N\right) }^{2^{j+1}\left( n+N\right)
-1}\left\vert b_{k}-b_{k+r}\right\vert +\sum\limits_{k=2^{j}\left(
n+N\right) }^{2^{j}\left( n+N\right) +r-1}\left\vert b_{k}\right\vert
\right\}
\end{equation*}%
\begin{equation*}
\ll \frac{N+1}{r}\sum\limits_{j=0}^{\infty }\left\{
\sum\limits_{k=2^{j}\left( n+N\right) }^{2^{j}\left( n+N\right)
+r-1}\left\vert b_{k}\right\vert +\sum\limits_{k=\left[ 2^{j}\left(
n+N\right) /c\right] }^{\left[ c2^{j}\left( n+N\right) \right] }\frac{%
\left\vert b_{k}\right\vert }{k}\right\}
\end{equation*}%
\begin{equation}
\ll \frac{N+1}{N+n}\varepsilon _{n}^{\left( 1\right)
}\sum\limits_{j=0}^{\infty }\frac{1}{2^{j}}\ll \varepsilon _{n}^{\left(
1\right) }.  \label{p12}
\end{equation}

Finally, we prove that (\ref{p1}) is true for $\frac{2l\pi }{r}+\frac{\pi }{r%
}\leq x<\frac{2\left( l+1\right) \pi }{r}$, where $0<2\left( l+1\right) \leq
r.$

Let $M:=M\left( x\right) \geq r$ be the natural number such that%
\begin{equation}
\frac{2\left( l+1\right) \pi }{r}-\frac{\pi }{M}\leq x<\frac{2\left(
l+1\right) \pi }{r}-\frac{\pi }{M+1}.  \label{p8}
\end{equation}%
Then%
\begin{equation*}
r_{n}\left( x\right) =\sum\limits_{k=n}^{n+M-1}b_{k}\sin
kx+\sum\limits_{k=n+M}^{\infty }b_{k}\sin kx=r_{n}^{\left( 3\right) }\left(
x\right) +r_{n}^{\left( 4\right) }\left( x\right) .
\end{equation*}%
Applying Lagrange's mean value theorem to the function $f\left( x\right)
=\sin kx$ on the interval $\left[ x,\frac{2\left( l+1\right) \pi }{r}\right] 
$ we obtain that there exists $z\in \left( x,\frac{2\left( l+1\right) \pi }{r%
}\right) $ such that%
\begin{equation*}
\sin \left( k\frac{2\left( l+1\right) \pi }{r}\right) -\sin kx=k\cos
kz\left( \frac{2\left( l+1\right) \pi }{r}-x\right) .
\end{equation*}%
Using this we get%
\begin{equation*}
r_{n}^{\left( 3\right) }\left( x\right)
=\sum\limits_{k=n}^{n+M-1}kb_{k}\cos kz\left( x-\frac{2\left( l+1\right)
\pi }{r}\right) +\sum\limits_{k=n}^{n+M-1}b_{k}\sin \left( k\frac{2\left(
l+1\right) \pi }{r}\right)
\end{equation*}%
\begin{equation*}
=r_{n}^{\left( 3.1\right) }\left( x\right) +r_{n}^{\left( 3.2\right) }\left( 
\frac{2\left( l+1\right) \pi }{r}\right) .
\end{equation*}%
Hence, by (\ref{p8}),%
\begin{equation}
\left\vert r_{n}^{\left( 3.1\right) }\left( x\right) \right\vert \leq \left( 
\frac{2\left( l+1\right) \pi }{r}-x\right)
\sum\limits_{k=n}^{n+M-1}k\left\vert b_{k}\right\vert \leq \pi \varepsilon
_{n}^{\left( 1\right) }.  \label{p14}
\end{equation}%
The quantity $r_{n}^{\left( 3.2\right) }\left( \frac{2\left( l+1\right) \pi 
}{r}\right) $ we can estimate in the same way as the quantity $r_{n}^{\left(
1.2\right) }\left( \frac{2l\pi }{r}\right) .$ Therefore we get%
\begin{equation}
\left\vert r_{n}^{\left( 3.2\right) }\left( \frac{2\left( l+1\right) \pi }{r}%
\right) \right\vert \ll \varepsilon _{n}^{\left( 1\right) }+\varepsilon
_{n}^{\left( 2\right) }.  \label{p15}
\end{equation}%
If $\left( b_{n}\right) \in GM\left( \beta ^{\ast },3\right) $ $\left( r\geq
3\right) $, then using Lemma 1, the inequality $2\left( l+1\right) -\frac{r}{%
\pi }x\leq \left\vert \sin \frac{rx}{2}\right\vert $ $\left( x\in \left[ 
\frac{2l\pi }{r}+\frac{\pi }{r},\frac{2\left( l+1\right) \pi }{r}\right] 
\text{ and }0<2\left( l+1\right) \leq r\right) $ and (\ref{p8}) we obtain%
\begin{equation*}
\left\vert r_{n}^{\left( 4\right) }\left( x\right) \right\vert =\left\vert
\sum\limits_{j=0}^{\infty }\sum\limits_{k=2^{j}\left( n+M\right)
}^{2^{j+1}\left( n+M\right) -1}b_{k}\sin kx\right\vert
\end{equation*}%
\begin{equation*}
\leq \sum\limits_{j=0}^{\infty }\left\vert \frac{-1}{2\sin \left(
rx/2\right) }\sum\limits_{k=2^{j}\left( n+M\right) }^{2^{j+1}\left(
n+M\right) -1}\left( b_{k}-b_{k+r}\right) \cos \left( k+\frac{r}{2}\right)
x\right.
\end{equation*}%
\begin{equation*}
\left. \left. +\sum\limits_{k=2^{j+1}\left( n+M\right) }^{2^{j+1}\left(
n+M\right) +r-1}b_{k}\cos \left( k-\frac{r}{2}\right)
x-\sum\limits_{k=2^{j}\left( n+M\right) }^{2^{j}\left( n+M\right)
+r-1}b_{k}\cos \left( k-\frac{r}{2}\right) x\right\} \right\vert
\end{equation*}%
\begin{equation*}
\leq \frac{1}{2\left\vert \sin \left( rx/2\right) \right\vert }%
\sum\limits_{j=0}^{\infty }\left\{ \sum\limits_{k=2^{j}\left( n+M\right)
}^{2^{j+1}\left( n+M\right) -1}\left\vert b_{k}-b_{k+r}\right\vert
+\sum\limits_{k=2^{j+1}\left( n+M\right) }^{2^{j+1}\left( n+M\right)
+r-1}\left\vert b_{k}\right\vert +\sum\limits_{k=2^{j}\left( n+M\right)
}^{2^{j}\left( n+M\right) +r-1}\left\vert b_{k}\right\vert \right\}
\end{equation*}%
\begin{equation*}
\ll \frac{1}{2\left( l+1\right) -\frac{r}{\pi }x}\sum\limits_{j=0}^{\infty
}\left\{ \sum\limits_{k=2^{j}\left( n+M\right) }^{2^{j+1}\left( n+M\right)
-1}\left\vert b_{k}-b_{k+r}\right\vert +\sum\limits_{k=2^{j}\left(
n+M\right) }^{2^{j}\left( n+M\right) +r-1}b_{k}\right\}
\end{equation*}%
\begin{equation*}
\ll \frac{M+1}{r}\sum\limits_{j=0}^{\infty }\left\{
\sum\limits_{k=2^{j}\left( n+M\right) }^{2^{j}\left( n+M\right)
+r-1}\left\vert b_{k}\right\vert +\sum\limits_{k=\left[ 2^{j}\left(
n+M\right) /c\right] }^{\left[ c2^{j}\left( n+M\right) \right] }\frac{%
\left\vert b_{k}\right\vert }{k}\right\}
\end{equation*}%
\begin{equation}
\ll \frac{M+1}{M+n}\varepsilon _{n}^{\left( 1\right)
}\sum\limits_{j=0}^{\infty }\frac{1}{2^{j}}\ll \varepsilon _{n}^{\left(
1\right) }.  \label{p16}
\end{equation}%
From the estimation (\ref{p9}), (\ref{p10}), (\ref{p11}), (\ref{p12}), (\ref%
{p14}), (\ref{p15}) and (\ref{p16}) we obtain the uniform convergence of
series (\ref{1}) follows and thus the proof is complete. $\square $

\subsection{Proof of Remark 4}

Let $r\geq 3$ and%
\begin{equation*}
a_{n}=\left\{ 
\begin{array}{c}
0\text{ \ \ \ \ \ \ \ \ \ \ \ if \ \ }r\mid n, \\ 
\frac{1}{n\ln \left( n+1\right) }\text{ \ \ if \ \ }r\nmid n.%
\end{array}%
\right.
\end{equation*}%
It is clear that $na_{n}\rightarrow 0$ as $n\rightarrow \infty .$

First, we prove that $\left( a_{n}\right) \in GM\left( \beta ^{\ast
},r\right) .$ Let%
\begin{equation*}
A_{r}:=A\left( r,k,n\right) =\left\{ k:n\leq k<2n\text{ and }r\mid k\right\}
\end{equation*}%
and%
\begin{equation*}
B_{r}:=B\left( r,k,n\right) \left\{ k:n\leq k<2n\text{ and }r\nmid k\right\}
.
\end{equation*}%
Then for all $n$%
\begin{equation*}
\sum\limits_{k=n}^{2n-1}\left\vert a_{k}-a_{k+r}\right\vert =\left(
\sum\limits_{k\in A_{r}}+\sum\limits_{k\in B_{r}}\right) \left\vert
a_{k}-a_{k+r}\right\vert
\end{equation*}%
\begin{equation*}
=\sum\limits_{k\in B_{r}}\left\vert \frac{1}{k\ln \left( k+1\right) }-\frac{%
1}{\left( k+r\right) \ln \left( k+r+1\right) }\right\vert
\end{equation*}%
\begin{equation*}
=\sum\limits_{k\in B_{r}}\frac{k\ln \left( 1+\frac{r}{k+1}\right) +r\ln
\left( k+r+1\right) }{k\left( k+r\right) \ln \left( k+1\right) \ln \left(
k+r+1\right) }\leq 2r\sum\limits_{k\in B_{r}}\frac{1}{k^{2}\ln \left(
k+1\right) }
\end{equation*}%
\begin{equation*}
=2r\sum\limits_{k=n}^{2n-1}\frac{a_{k}}{k}\ll \sum\limits_{k=\left[ n/c%
\right] }^{\left[ cn\right] }\frac{a_{k}}{k}
\end{equation*}%
and $\left( a_{n}\right) \in GM\left( \beta ^{\ast },r\right) $.

Now, we show that (\ref{r1}) is valid. We have%
\begin{equation*}
\sum\limits_{n=1}^{\infty }\sum\limits_{k=1}^{\left[ r/2\right]
}\left\vert a_{r\cdot n+k}-a_{r\cdot n+r-k}\right\vert
\end{equation*}%
\begin{equation*}
=\sum\limits_{n=1}^{\infty }\sum\limits_{k=1}^{\left[ r/2\right]
}\left\vert \frac{1}{\left( rn+k\right) \ln \left( rn+k+1\right) }-\frac{1}{%
\left( rn+r-k\right) \ln \left( rn+r-k+1\right) }\right\vert
\end{equation*}%
\begin{equation*}
\leq \sum\limits_{n=1}^{\infty }\sum\limits_{k=1}^{\left[ r/2\right] \ }%
\frac{rn\ln \left( 1+\frac{r-2k}{rn+k+1}\right) +r\ln \left( rn+r-k+1\right) 
}{\left( rn+k\right) \left( rn+r-k\right) \ln \left( rn+k+1\right) \ln
\left( rn+r-k+1\right) }
\end{equation*}%
\begin{equation*}
\leq \sum\limits_{n=1}^{\infty }\sum\limits_{k=1}^{\left[ r/2\right] }%
\frac{r-2k+r\ln \left( rn+r-k+1\right) }{\left( rn+k\right) \left(
rn+r-k\right) \ln \left( rn+k+1\right) \ln \left( rn+r-k+1\right) }
\end{equation*}%
\begin{equation*}
\leq 2r\sum\limits_{n=1}^{\infty }\sum\limits_{k=1}^{\left[ r/2\right] }%
\frac{1}{\left( rn+k\right) \left( rn+r-k\right) \ln \left( rn+k+1\right) }%
\leq \sum\limits_{n=1}^{\infty }\frac{1}{n^{2}\ln \left( n+1\right) }
\end{equation*}%
and since the series $\sum\limits_{n=1}^{\infty }\frac{1}{n^{2}\ln \left(
n+1\right) }$ converges (\ref{r1}) holds.

Finally, we prove that $\left( a_{n}\right) \notin GM\left( \beta ^{\ast
},2\right) $. For all $n$ we get%
\begin{equation*}
\sum\limits_{k=n}^{2n-1}\left\vert a_{k}-a_{k+2}\right\vert \geq
\sum\limits_{k\in A_{r}}\left\vert a_{k}-a_{k+2}\right\vert
=\sum\limits_{k\in A_{r}}\frac{1}{\left( k+2\right) \ln \left( k+3\right) }
\end{equation*}%
\begin{equation*}
\geq \frac{1}{3r\ln \left( 2n+2\right) }
\end{equation*}%
and since%
\begin{equation*}
\sum\limits_{k=n}^{n+1}a_{k}+\sum\limits_{k=\left[ n/c\right] }^{\left[ cn%
\right] }\frac{a_{k}}{k}\ll \frac{1}{n\ln \left( n+1\right) },
\end{equation*}%
the inequality%
\begin{equation*}
\sum\limits_{k=n}^{2n-1}\left\vert a_{k}-a_{k+2}\right\vert \leq K\left(
a\right) \left\{ \sum\limits_{k=n}^{n+1}a_{k}+\sum\limits_{k=\left[ n/c%
\right] }^{\left[ cn\right] }\frac{a_{k}}{k}\right\}
\end{equation*}%
does not hold, that is, $\left( a_{n}\right) $ does not belong to $GM\left(
\beta ^{\ast },2\right) .$

This ends our proof. $\square $

\end{document}